\documentclass[12pt,a4paper]{article}
\usepackage[utf8x]{inputenc}
\usepackage{amsmath}
\usepackage{amsfonts}
\usepackage{amssymb}
\usepackage{amsthm}
\usepackage{mathrsfs}
\usepackage{fancyhdr}
\usepackage{graphicx}
\usepackage[usenames,x11names]{xcolor}
\usepackage{arabtex}
\usepackage{framed}
\usepackage[top=2cm,bottom=2cm,margin=2cm]{geometry}
\usepackage{cancel}
\usepackage{array}   
\usepackage{ulem}

\usepackage{hyperref}

\title{Characterization of the Bernoulli polynomials via the Raabe functional equation}
\author{\sc Bakir FARHI \\
Laboratoire de Mathématiques appliquées \\
Faculté des Sciences Exactes \\
Université de Bejaia, 06000 Bejaia, Algeria \\[1mm]
\href{mailto:bakir.farhi@univ-bejaia.dz}{bakir.farhi@univ-bejaia.dz} \\[1mm]
\url{http://farhi.bakir.free.fr/}
}
\date{}

\let\up=\textsuperscript

\def\R{{\mathbb R}}
\def\Q{{\mathbb Q}}
\def\C{{\mathbb C}}
\def\N{{\mathbb N}}
\def\Z{{\mathbb Z}}
\def\E{{\mathbb E}}

\def\odd{\mathrm{odd}} 

\def\ind{1\hspace*{-0.9mm}\mathrm{I}} 

\def\restmod#1#2{#1\ (\mathrm{mod}\ #2)} 

\def\idem{\leavevmode\hbox to 10.6mm{\vrule height .63ex depth -.59ex
    width 10mm\hfill}}

\theoremstyle{plain}
\numberwithin{equation}{section}
\newtheorem{thm}{Theorem}[section]

\newtheorem{lemma}[thm]{Lemma}

\newtheorem{coll}[thm]{Corollary}

\theoremstyle{definition}

\theoremstyle{remark}
\newtheorem{rmk}[thm]{Remark}
\newtheorem{rmks}[thm]{Remarks}
\newtheorem{expl}[thm]{Example}
\newtheorem{expls}[thm]{Examples}

\begin{document}
\maketitle

\begin{abstract}
The purpose of the present paper is to show that in certain classes of real (or complex) functions, the Bernoulli polynomials are essentially the only ones satisfying the Raabe functional equation. For the class of the real $1$-periodic functions which are expandable as Fourier series, we point out new solutions of the Raabe functional equation, not relating to the Bernoulli polynomials. Furthermore, we will give for the considered classes various proofs, making the mathematical content of the paper quite rich.
\end{abstract}

\noindent\textbf{MSC 2020:} Primary 11B68; Secondary 39B22, 39B32, 30D05. \\
\textbf{Keywords:} The Raabe multiplicative theorem, Bernoulli polynomials, functional equations for real or complex functions.

\section{Introduction and Notation}\label{sec1}

Throughout this paper, we let $\N$ denote the set of positive integers and $\N_0 := \N \cup \{0\}$ the set of nonnegative integers. For $x \in \R$, we let $\lfloor{x}\rfloor$ and $\langle{x}\rangle$ respectively denote the integer part and the fractional part of $x$. Next, we let $B_n(X)$ ($n \in \N_0$) denote \textit{the Bernoulli polynomials} which can be defined through their exponential generating function:
\begin{equation}\label{eq1}
\frac{t e^{X t}}{e^t - 1} = \sum_{n = 0}^{+ \infty} B_n(X) \frac{t^n}{n!} .
\end{equation} 
The Bernoulli polynomials, together with their values at $X = 0$ (called \textit{the Bernoulli numbers}), have many important properties and appear in several areas of mathematics, including number theory, real analysis, special functions, combinatorics, orthogonal polynomials, and so on. We refer to the book of Nielsen \cite{nie} for an elementary (but old) presentation and to the lecture notes of Kouba \cite{kou} for a modern presentation. It is known for example that for every $n \in \N_0$ the Bernoulli polynomial $B_n(X)$ is monic of degree $n$ and satisfies the following polynomial identities ((see \cite{kou})):
\begin{align}
B_n'(X) & = n B_{n - 1}(X) , \label{eq4} \\
B_n(X + 1) - B_n(X) & = n X^{n - 1} . \label{eq14}
\end{align}
It is also known that the Fourier series expansion of the Bernoulli polynomial $B_n(X)$ ($n \in \N$) in the interval $[0 , 1)$ is given by (see \cite{kou}):
\begin{equation}\label{eq18}
B_n(x) = \frac{- 2 \cdot n!}{(2 \pi)^n} \sum_{k = 1}^{+ \infty} \frac{\cos\left(2 \pi k x - \frac{n \pi}{2}\right)}{k^n} 
\end{equation}
(which is valid for all $x \in [0 , 1)$ if $n \geq 2$ and all $x \in (0 , 1)$ if $n = 1$). By analogy with \eqref{eq18}, let us define
\begin{equation}\label{eq19}
\overline{B}_n(x) = \frac{- 2 \cdot n!}{(2 \pi)^n} \sum_{k = 1}^{+ \infty} \frac{\sin\left(2 \pi k x - \frac{n \pi}{2}\right)}{k^n} 
\end{equation}
(for all $(n , x) \in \N \times [0 , 1) \setminus \{(1 , 0)\}$).

In \cite{raa}, Raabe showed that for every $n \in \N_0$ and every $a \in \N$, the Bernoulli polynomial $B_n(X)$ satisfies the remarkable functional equation (in $f$):
\begin{equation}
f(a x) = a^{n - 1} \left(f(x) + f\left(x + \frac{1}{a}\right) + f\left(x + \frac{2}{a}\right) + \dots + f\left(x + \frac{a - 1}{a}\right)\right) . \tag{$\E_{n , a}$}
\end{equation} 
(Remark that this equation is trivial for $a = 1$). Nielsen \cite{nie} remarked that for every fixed integer $a \geq 2$, the Bernoulli polynomial $B_n(X)$ is the only monic polynomial satisfying $(\E_{n , a})$; a simple proof of this result will be given later in this paper. Carlitz \cite{car2} generalized the Raabe theorem by showing that for all $n \in \N_0$ and all $a , b \in \N$, we have:
\begin{equation}\label{eq2}
a^{n - 1} \sum_{k = 0}^{a - 1} B_n\left(\frac{x}{a} + \frac{b k}{a}\right) = b^{n - 1} \sum_{\ell = 0}^{b - 1} B_n\left(\frac{x}{b} + \frac{a \ell}{b}\right) .
\end{equation}
Notice that \eqref{eq2} is symmetric in $a$ and $b$ and that the Raabe theorem can be derived from it by taking $b = 1$ and replacing $x$ by $a x$. In \cite{car3}, Carlitz showed in addition that for every $n \in \N_0$, the unique monic polynomial satisfying \eqref{eq2} for two different values $a$ and $b$ is the Bernoulli polynomial $B_n(X)$. Related results concerning not only the Bernoulli polynomials but also the Euler polynomials are given in \cite{car1,car2,car3}.    

In This paper, given $n \in \N_0$ and $a \geq 2$ an integer, we show that the Bernoulli polynomial $B_n(X)$ is essentially the only solution (up to a multiplicative constant) of Equation $(\E_{n , a})$ in certain classes of real (or complex) functions. The considered classes of functions are, in this order, the complex polynomials, the complex analytic functions, $\mathscr{C}^n(\R , \R)$, $\mathscr{C}^0(\R , \R)$, and the real $1$-periodic functions which are expandable as Fourier series. Exceptionally for the last class, the set of solutions of $(\E_{n , a})$ contains new solutions not (directly) related to the Bernoulli polynomial $B_n(X)$ (see Theorem \ref{t5} and Corollary \ref{coll2}). Further, the methods used to solve Equation $(\E_{n , a})$ in each of the considered classes of functions are various, which makes the mathematical content of the paper quite rich.

\section{Some useful lemmas about Equation $(\E_{n , a})$}\label{sec2}

The first lemma below is obvious and needs no proof.

\begin{lemma}\label{l1}
Let $n \in \N_0$ and $a \in \N$. Then the set of functions $f : \R \rightarrow \R$ which are solutions of Equation $(\E_{n , a})$ constitutes an $\R$-vector space. Similarly, the set of functions $f : \C \rightarrow \C$ which are solutions of Equation $(\E_{n , a})$ constitutes a $\C$-vector space. \hfill $\square$ 
\end{lemma}

The next two lemmas are crucial for proving our main results (see \S\ref{sec3}); they provide important properties for the spaces of solutions of Equations $(\E_{n , a})$.

\begin{lemma}\label{l2}
Let $n , a \in \N$. If a differentiable real {\rm(}or complex{\rm)} function $f$ is a solution of Equation $(\E_{n , a})$ then its derivative $f'$ is a solution of Equation $(\E_{n - 1 , a})$. 
\end{lemma}

\begin{proof}
Let $f$ be a differentiable real (or complex) function. Suppose that $f$ is a solution of Equation $(\E_{n , a})$; that is
$$
f(a x) = a^{n - 1} \left[f(x) + f\left(x + \frac{1}{a}\right) + f\left(x + \frac{2}{a}\right) + \dots + f\left(x + \frac{a - 1}{a}\right)\right] .
$$
By differentiating this last identity, we get (after simplifying and rearranging):
$$
f'(a x) = a^{n - 2} \left[f'(x) + f'\left(x + \frac{1}{a}\right) + f'\left(x + \frac{2}{a}\right) + \dots + f'\left(x + \frac{a - 1}{a}\right)\right] ,
$$
implying that $f'$ is a solution of Equation $(\E_{n - 1 , a})$, as required.
\end{proof}

\begin{lemma}\label{l3}
Let $n \in \N_0$ and $a , b\in \N$. If a real {\rm(}or complex{\rm)} function $f$ is both a solution of Equation $(\E_{n , a})$ and Equation $(\E_{n , b})$ then it is a solution of Equation $(\E_{n , a b})$.
\end{lemma}

\begin{proof}
Suppose that a real (resp. complex) function $f$ is both a solution of Equation $(\E_{n , a})$ and Equation $(\E_{n , b})$. Then we have for all $x \in \R$ (resp. $x \in \C$):
\begin{equation}\label{eq3}
f(a b x) = f(a \cdot b x) = a^{n - 1} \sum_{k = 0}^{a - 1} f\left(b x + \frac{k}{a}\right) .
\end{equation}
Next, since $f$ satisfies Equation $(\E_{n , b})$ then we have for all $x \in \R$ (resp. $x \in \C$) and all $k \in \{0 , 1 , \dots , a - 1\}$:
\begin{align*}
f\left(b x + \frac{k}{a}\right) & = f\left(b \cdot \left(x + \frac{k}{a b}\right)\right) \\
& = b^{n - 1} \sum_{\ell = 0}^{b - 1} f\left(x + \frac{k}{a b} + \frac{\ell}{b}\right) \\
& = b^{n - 1} \sum_{\ell = 0}^{b - 1} f\left(x + \frac{a \ell + k}{a b}\right) .
\end{align*}
By inserting this (for $0 \leq k \leq a - 1$) into \eqref{eq3}, we get:
\begin{align*}
f(a b x) & = a^{n - 1} \sum_{k = 0}^{a - 1} b^{n - 1} \sum_{\ell = 0}^{b - 1} f\left(x + \frac{a \ell + k}{a b}\right) \\
& = (a b)^{n - 1} \sum_{\begin{subarray}{c}
0 \leq k \leq a - 1 \\
0 \leq \ell \leq b - 1
\end{subarray}} f\left(x + \frac{a \ell + k}{a b}\right) .
\end{align*}
But since the sequence of integers ${(a \ell + k)}_{\begin{subarray}{c}
0 \leq k \leq a - 1 \\
0 \leq \ell \leq b - 1
\end{subarray}}$ is nothing else the sequence $0 , 1 , \dots , a b - 1$ (up to an order), we conclude to
$$
f(a b x) = (a b)^{n - 1} \sum_{r = 0}^{a b - 1} f\left(x + \frac{r}{a b}\right) ,
$$
for all $x \in \R$ (resp. $x \in \C$). This shows that $f$ is a solution of Equation $(\E_{n , a b})$, as required.
\end{proof}

\begin{rmk}
Let $n \in \N_0$ be fixed. According to Lemma \ref{l3}, for a function $f$ to be a solution of all Equations $(\E_{n , a})$ ($a \in \N$), it suffices that $f$ is a solution of all Equations $(\E_{n , p})$ ($p$ prime).
\end{rmk}

\section{The main results}\label{sec3}

The first result below is already pointed out by Nielsen \cite{nie}. Here, we give for it an alternative nice proof.

\begin{thm}\label{t1}
Let $n \in \N_0$ and $a \geq 2$ an integer. Then the solutions of Equation $(\E_{n , a})$ in $\C[X]$ are the polynomials $\lambda B_n(X)$ {\rm(}$\lambda \in \C${\rm)}. 
\end{thm}

\begin{proof}
By the Raabe multiplication theorem, we know that the polynomials $\lambda B_n(X)$ ($\lambda \in \C$) are all solutions of Equation $(\E_{n , a})$. To prove the converse, we argue by induction on $n \in \N_0$. \\[1mm]
\textbullet{} For $n = 0$. Let $P \in \C[X]$ be a solution of Equation $(\E_{0 , a})$; that is
$$
P(a X) = \frac{1}{a} \left[P(X) + P\left(X + \frac{1}{a}\right) + P\left(X + \frac{2}{a}\right) + \dots + P\left(X + \frac{a - 1}{a}\right)\right] .
$$
If $P = 0_{\C[X]}$ then we are done (since $0_{\C[X]} = 0 B_0(X)$). Suppose now that $P \neq 0_{\C[X]}$ and let $d \in \N_0$ denote the degree of $P$ and $\lambda \in \C^*$ denote the leading coefficient of $P$. By comparing the leading coefficients of the two sides of the last polynomial identity, we get
$$
\lambda a^d = \lambda , 
$$
giving $d = 0$. This means that $P$ is constant, that is $P(X) = \lambda = \lambda B_0(X)$, as required. \\[1mm]
\textbullet{} Let $n \in \N$. Suppose that every solution in $\C[X]$ of Equation $(\E_{n - 1 , a})$ has the form $\lambda B_{n - 1}(X)$ ($\lambda \in \C$) and let us show that every solution in $\C[X]$ of Equation $(\E_{n , a})$ has the form $\lambda B_n(X)$ ($\lambda \in \C$). Let $P \in \C[X]$ be a solution of Equation $(\E_{n , a})$; so (according to Lemma \ref{l2}) its derivative $P' \in \C[X]$ is a solution of Equation $(\E_{n - 1 , a})$. Thus, according to our induction hypothesis, there exists $\lambda \in \C$ such that $P'(X) = \lambda B_{n - 1}(X)$. It follows by integrating and taking into account \eqref{eq4} that: $P(X) = \frac{\lambda}{n} B_n(X) + c$ (for some $c \in \C$). So, it remains to show that $c = 0$. Since $P(X)$ and $B_n(X)$ are both solutions of Equation $(\E_{n , a})$ then their linear combination $\left(P(X) - \frac{\lambda}{n} B_n(X)\right) = c$ is also a solution of $(\E_{n , a})$ (according to Lemma \ref{l1}), meaning that:
$$
c = a^{n - 1} \left(\underbrace{c + c + \dots + c}_{a \text{ times}}\right) = a^n c ,
$$
implying that $c = 0$, as required. This completes this induction and achieves the proof.
\end{proof}

\begin{thm}\label{t2}
Let $n \in \N_0$ and $a \geq 2$ be an integer. Then the only complex analytic functions on $\C$ which are solutions of Equation $(\E_{n , a})$ are the polynomials $\lambda B_n(X)$ {\rm(}$\lambda \in \C${\rm)}.
\end{thm}

\begin{proof}
Since the polynomials $\lambda B_n(X)$ ($\lambda \in \C$) are all solutions of Equation $(\E_{n , a})$ (by the Raabe theorem) then we have just to show that any complex analytic function on $\C$ which satisfies $(\E_{n , a})$ has the form $\lambda B_n(X)$ ($\lambda \in \C$). We shall show this fact in two steps. \\[1mm]
\textbullet{} \textbf{1\up{st} step.} (The case $n = 0$). Let $f$ be a complex analytic function on $\C$ which satisfies Equation $(\E_{0 , a})$, that is
\begin{equation}\label{eq5}
f(a z) = \frac{1}{a} \left[f(z) + f\left(z + \frac{1}{a}\right) + f\left(z + \frac{2}{a}\right) + \dots + f\left(z + \frac{a - 1}{a}\right)\right] ~~~~~~~~~~ (\forall z \in \C) ,
\end{equation}
and let us show that $f$ is necessarily constant. Since $f$ is analytic (so continuous) on $\C$ then it is in particular continuous on the closed ball $\overline{B}(0 , 2) := \{u \in \C : |u| \leq 2\}$ of $\C$. Next, since $\overline{B}(0 , 2)$ is compact (it is closed and bounded) then $f$ is bounded on $\overline{B}(0 , 2)$. So, let us define
$$
M := \sup_{\begin{subarray}{c}
u \in \C \\
|u| \leq 2
\end{subarray}} \left\vert{f(u)}\right\vert \in [0 , + \infty) .
$$
By leaning on \eqref{eq5}, we are going to show that $f$ is even bounded on $\C$ by $M$, which implies (according to Liouville's theorem) that $f$ is constant. Let $z \in \C$ arbitrary. Consider a positive integer $k$ such that $a^k \geq |z|$ (such a $k$ exists since $\lim_{m \rightarrow + \infty} a^m = + \infty$) and set $z_0 := \frac{z}{a^k}$. So we have $|z_0| \leq 1$. Next, since $f$ satisfies Equation $(\E_{0 , a})$ then (by iterating the result of Lemma \ref{l3}) it satisfies all Equations $(\E_{0 , a^r})$ ($r \in \N$); in particular $f$ satisfies Equation $(\E_{0 , a^k})$. Thus we have
$$
f(z) = f\left(a^k z_0\right) = \frac{1}{a^k} \sum_{\ell = 0}^{a^k - 1} f\left(z_0 + \frac{\ell}{a^k}\right) .
$$
It follows by using the triangle inequality in $\C$ that:
$$
\left\vert{f(z)}\right\vert \leq \frac{1}{a^k} \sum_{0 \leq \ell \leq a^k - 1} \left\vert{f\left(z_0 + \frac{\ell}{a^k}\right)}\right\vert .
$$
But for every $\ell \in \{0 , 1 , \dots , a^k - 1\}$, we have that $\left\vert{z_0 + \frac{\ell}{a^k}}\right\vert \leq \left\vert{z_0}\right\vert + \frac{\ell}{a^k} \leq 1 + 1 = 2$, implying that $\left\vert{f\left(z_0 + \frac{\ell}{a^k}\right)}\right\vert \leq M$. It follows from this fact that:
$$
\left\vert{f(z)}\right\vert \leq \frac{1}{a^k} \sum_{0 \leq \ell \leq a^k - 1} M = M .
$$
Consequently, $f$ is bounded on $\C$; thus (by Liouville's theorem) it is constant, as required. \\[1mm]
\textbullet{} \textbf{2\up{nd} step.} (The general case). Let $f$ be a complex analytic function on $\C$ which satisfies Equation $(\E_{n , a})$. By iterating the result of Lemma \ref{l2}, we derive that the $n$\up{th} derivative $f^{(n)}$ of $f$ is a solution of Equation $(\E_{0 , a})$, implying (according to the result of the first step above) that $f^{(n)}$ is constant. Thus $f$ is a polynomial of degree $\leq n$, which implies (according to Theorem \ref{t1}) that $f$ has the required form $f(X) = \lambda B_n(X)$ ($\lambda \in \C$). This completes the proof of the theorem. 
\end{proof}

\begin{thm}\label{t3}
Let $a \geq 2$ be an integer. Then the only solutions of Equation $(\E_{0 , a})$ in $\mathscr{C}^0(\R , \R)$ are the constant functions.
\end{thm}

\begin{proof}
It is obvious that the constant functions are solutions of Equation $(\E_{0 , a})$. Conversely, let $f \in \mathscr{C}^0(\R , \R)$ be a solution of $(\E_{0 , a})$ and let us prove that $f$ is constant. According to Lemma \ref{l3}, $f$ is a solution of all Equations $(\E_{0 , a^k})$ ($k \in \N$); that is for all $k \in \N$ and all $x \in \R$, we have
\begin{equation}\label{eq6}
f\left(a^k x\right) = \frac{1}{a^k} \left[f(x) + f\left(x + \frac{1}{a^k}\right) + f\left(x + \frac{2}{a^k}\right) + \dots + f\left(x + \frac{a^k - 1}{a^k}\right)\right] .
\end{equation}
We remark that the right-hand side of \eqref{eq6} is nothing else the left Riemann sum of $f$ with respect to the regular subdivision $x_0 = x$, $x_1 = x + \frac{1}{a^k}$, $x_2 = x + \frac{2}{a^k}$, \dots, $x_{a^k} = x + 1$ of the interval $[x , x + 1]$. Since $f \in \mathscr{C}^0(\R , \R)$ then it is Riemann integrable on every bounded closed interval of $\R$; in particular, $f$ is Riemann integrable on $[x , x + 1]$ for every $x \in \R$. So we derive from \eqref{eq6} that we have for all $x \in \R$:
\begin{equation}\label{eq7}
\lim_{k \rightarrow + \infty} f\left(a^k x\right) = \int_{x}^{x + 1} f(t) \, d t .
\end{equation}
Now, by using \eqref{eq7}, we have for all $x \in \R$ and all $\ell \in \N$:
\begin{align*}
\int_{x}^{x + 1} f(t) \, d t & = \lim_{k \rightarrow + \infty} f\left(a^k x\right) \\
& = \lim_{k \rightarrow + \infty} f\left(a^{k - \ell} x\right) \\
& = \lim_{k \rightarrow + \infty} f\left(a^k \left(\frac{x}{a^{\ell}}\right)\right) \\
& = \int_{\frac{x}{a^{\ell}}}^{\frac{x}{a^{\ell}} + 1} f(t) \, d t ~~~~~~~~~~ \left(\text{by applying \eqref{eq7} for } \frac{x}{a^{\ell}} \text{ instead of } x\right) .
\end{align*}
By tending $\ell$ to infinity and taking into account the continuity of $f$, it follows that we have for all $x \in \R$:
\begin{equation}\label{eq8}
\int_{x}^{x + 1} f(t) \, d t = \int_{0}^{1} f(t) \, d t .
\end{equation}
Then by differentiating, we get
$$
f(x + 1) - f(x) = 0 ~~~~~~~~~~ (\forall x \in \R) ,
$$
meaning that $f$ is $1$-periodic. Next, set
\begin{align*}
\sigma & := \int_{0}^{1} f(t) \, d t \\[-7mm]
\intertext{and} \\[-14mm]
E & := \left\{\frac{r}{a^s - 1} ;~ r \in \Z , s \in \N\right\} \subset \R .
\end{align*}
Our strategy now consists in showing that $f$ is constant on $E$ (equal to $\sigma$) and that $E$ is dense in $\R$, concluding that $f$ is constant on $\R$ (since it is continuous on $\R$). According to \eqref{eq7} and \eqref{eq8}, we have for all $x \in \R$:
\begin{equation}\label{eq9}
\lim_{k \rightarrow + \infty} f\left(a^k x\right) = \sigma .
\end{equation}
Let us first show that $f$ is constant on $E$; more precisely that $f(x) = \sigma$, $\forall x \in E$. Let $x \in E$ arbitrary. So we can write $x = \frac{r}{a^s - 1}$ for some $r \in \Z$ and some $s \in \N$. We remark that for all $k \in \N$, we have $(a^{k s} - 1) x \in \Z$ (since the integer $((a^{k s} - 1)$ is a multiple of the integer $(a^s - 1)$); that is $a^{k s} x - x \in \Z$. But since $f$ is $1$-periodic (as proved above), we derive that
$$
f(x) = f\left(a^{k s} x\right) ~~~~~~~~~~ (\forall k \in \N) .
$$
It follows by tending $k$ to $+ \infty$ and using \eqref{eq9} that:
$$
f(x) = \sigma .
$$
Thus $f$ is constant (equal to $\sigma$) on $E$, as claimed it to be. Now, let us show that $E$ is dense in $\R$. So we have to show that for all $u \in \R$, there exists a sequence ${(u_n)}_{n \in \N}$ of $E$ which converges to $u$. Let $u \in \R$ arbitrary. For all $n \in \N$, define
$$
r_n := \left\lfloor\left(a^n - 1\right) u\right\rfloor ~~~~\text{and}~~~~ u_n := \frac{r_n}{a^n - 1} .
$$
It is clear that ${(u_n)}_n$ is a sequence of $E$. In addition, we have for all $n \in \N$:
\begin{align*}
u - u_n & = u - \frac{r_n}{a^n - 1} \\[1mm]
& = \frac{\left(a^n - 1\right) u - \left\lfloor\left(a^n - 1\right) u\right\rfloor}{a^n -1} \\[1mm]
& = \frac{\left\langle\left(a^n - 1\right) u\right\rangle}{a^n - 1} .
\end{align*}
Thus
$$
0 \leq u - u_n < \frac{1}{a^n - 1} ~~~~~~~~~~ (\forall n \in \N) ,
$$
implying that ${(u_n)}_n$ converges to $u$ (since $\frac{1}{a^n - 1} \rightarrow 0$ as $n \rightarrow + \infty$). Hence $E$ is dense in $\R$, as claimed it to be. Finally, the three facts: ``$f$ is continuous on $\R$'', ``$f$ is constant on $E$'', and ``$E$ is dense in $\R$'' conclude that $f$ is constant on $\R$, as required. The proof is complete.
\end{proof}

\begin{rmk}
There exist solutions of Equation $(\E_{0 , a})$ ($a \in \N$, $a \geq 2$) which are not constant. We may think, for example, of the real function $\ind_{\Q}$ (the indicator function of the set of rational numbers $\Q$).
\end{rmk}

We now derive from Theorem \ref{t3} the following more general result:

\begin{coll}\label{coll1}
Let $n \in \N_0$ and $a \geq 2$ be an integer. Then the only solutions of Equation $(\E_{n , a})$ in $\mathscr{C}^n(\R , \R)$ are the polynomials $\lambda B_n(X)$ {\rm(}$\lambda \in \R${\rm)}.
\end{coll}

\begin{proof}
Let $f \in \mathscr{C}^n(\R , \R)$ be a solution of Equation $(\E_{n , a})$. Then (by iterating the result of Lemma \ref{l2}) the $n$\up{th} derivative $f^{(n)}$ of $f$ (which belongs to $\mathscr{C}^0(\R, \R)$ since $f$ belongs to $\mathscr{C}^n(\R , \R)$) is a solution of Equation $(\E_{0 , a})$. It follows (according to Theorem \ref{t3}) that $f^{(n)}$ is constant. Thus $f$ is a polynomial, which implies (according to Theorem \ref{t1}) that $f$ has the required form $f(X) = \lambda B_n(X)$ ($\lambda \in \R$). This achieves the proof.
\end{proof}

We are now interested in the solutions of Equations $(\E_{n , a})$ in $\mathscr{C}^0(\R , \R)$. We have the following theorem:

\begin{thm}\label{t4}
Let $n \in \N_0$ and $a \geq 2$ be an integer. Let $f \in \mathscr{C}^0(\R , \R)$ such that the limit $\lim_{h \rightarrow 0} \frac{\int_{h}^{h + 1} f(t) \, d t}{h^n}$ exists {\rm(}possibly infinite{\rm)}. Then $f$ is a solution of Equation $(\E_{n , a})$ if and only if it has the form
\begin{equation}\label{eq10}
f(X) = \lambda B_n(X) + \tau(X) ,
\end{equation}
where $\lambda \in \R$ and $\tau$ is a continuous real $1$-periodic function which is a solution of $(\E_{n , a})$.
\end{thm}

\begin{proof}
It is immediate that the functions having the form in \eqref{eq10} are solutions of $(\E_{n , a})$. Conversely, let $f \in \mathscr{C}^0(\R , \R)$ be a solution of $(\E_{n , a})$ and let us show that it has the form in \eqref{eq10}. According to Lemma \ref{l3}, $f$ is a solution of all Equations $(\E_{n , a^k})$ ($k \in \N$); that is for all $k \in \N$ and all $x \in \R$, we have
$$
f\left(a^k x\right) = a^{k (n - 1)} \left[f(x) + f\left(x + \frac{1}{a^k}\right) + f\left(x + \frac{2}{a^k}\right) + \dots + f\left(x + \frac{a^k - 1}{a^k}\right)\right] ;
$$
or equivalently
\begin{equation}\label{eq11}
\frac{f\left(a^k x\right)}{a^{k n}} = \frac{1}{a^k} \sum_{0 \leq \ell < a^k} f\left(x + \frac{\ell}{a^k}\right) .
\end{equation}
We remark that the right-hand side of \eqref{eq11} is exactly the left Riemann sum of $f$ with respect to the regular subdivision $x_0 = x$, $x_1 = x + \frac{1}{a^k}$, $x_2 = x + \frac{2}{a^k}$, \dots, $x_{a^k} = x + 1$ of the interval $[x , x + 1]$. Since $f \in \mathscr{C}^0(\R , \R)$ then it is Riemann integrable on every bounded closed interval of $\R$; in particular, $f$ is Riemann integrable on $[x , x + 1]$ for every $x \in \R$. So by tending $k$ to $+ \infty$ in \eqref{eq11}, we derive that we have for all $x \in \R$:
\begin{equation}\label{eq12}
\lim_{k \rightarrow + \infty} \frac{f\left(a^k x\right)}{a^{k n}} = \int_{x}^{x + 1} f(t) \, d t .
\end{equation}
Now, by using \eqref{eq12}, we have for all $x \in \R^*$ and all $\ell \in \N$:
\begin{align*}
\int_{x}^{x + 1} f(t) \, d t & = \lim_{k \rightarrow + \infty} \frac{f\left(a^k x\right)}{a^{k n}} \\
& = \lim_{k \rightarrow + \infty} \frac{f\left(a^{k - \ell} x\right)}{a^{(k - \ell) n}} \\
& = a^{\ell n} \lim_{k \rightarrow + \infty} \frac{f\left(a^k\left(\frac{x}{a^{\ell}}\right)\right)}{a^{k n}}
\end{align*}

\begin{align*}
& = a^{\ell n} \int_{\frac{x}{a^{\ell}}}^{\frac{x}{a^{\ell}} + 1} f(t) \, d t ~~~~~~~~~~ (\text{by applying \eqref{eq12} for } \frac{x}{a^{\ell}} \text{ instead of } x) \\
& = x^n \left(\frac{a^{\ell}}{x}\right)^n \int_{\frac{x}{a^{\ell}}}^{\frac{x}{a^{\ell}} + 1} f(t) \, d t . 
\end{align*}
Then, by tending $\ell$ to infinity and setting $h = \frac{x}{a^{\ell}}$ in the last expression, we get for all $x \in \R^*$:
$$
\int_{x}^{x + 1} f(t) \, d t = x^n \lim_{h \rightarrow 0} \frac{\int_{h}^{h + 1} f(t) \, d t}{h^n} . 
$$
By specializing in particular $x$ to $1$ in this last equality, we obtain that
$$
\lim_{h \rightarrow 0} \frac{\int_{h}^{h + 1} f(t) \, d t}{h^n} = \int_{1}^{2} f(t) \, d t .
$$
So, setting
$$
\sigma := \int_{1}^{2} f(t) \, d t ,
$$
we derive that we have for all $x \in \R^*$:
\begin{equation}\label{eq13}
\int_{x}^{x + 1} f(t) \, d t = \sigma x^n .
\end{equation}
But in view of the continuity of $f$, \eqref{eq13} is still true for $x = 0$, so it is true for all $x \in \R$. Then, by differentiating \eqref{eq13}, we get for all $x \in \R$:
$$
f(x + 1) - f(x) = n \sigma x^{n - 1} ;
$$
that is (according to \eqref{eq14}):
$$
f(x + 1) - f(x) = \sigma \left(B_n(x + 1) - B_n(x)\right) .
$$
So, by defining a function $\tau$ as:
$$
\tau(x) := f(x) - \sigma B_n(x) ~~~~~~~~~~ (\forall x \in \R) , 
$$
we have that
$$
\tau(x + 1) = \tau(x) ~~~~~~~~~~ (\forall x \in \R) ,
$$
showing that $\tau$ is $1$-periodic. In addition, $\tau$ is continuous and it is a solution of Equation $(\E_{n , a})$ (since it is a linear combination of $f$ and $B_n$, which are both continuous and solutions of Equation $(\E_{n , a})$). Consequently, $f$ has the required form $f(X) = \sigma B_n(X) + \tau(X)$. This completes the proof.
\end{proof}

\begin{rmks}~
\begin{enumerate}
\item The condition of the existence of the limit $\lim_{h \rightarrow 0} \frac{\int_{h}^{h + 1} f(t) \, d t}{h^n}$ occurring in Theorem \ref{t4} is automatically verified for $n = 0$ (since $f$ is continuous). This is why this condition doesn't appear in Theorem \ref{t3}.  
\item An alternative condition which can replace the condition ``$\lim_{h \rightarrow 0} \frac{\int_{h}^{h + 1} f(t) \, d t}{h^n}$ exists'' in Theorem \ref{t4} is ``$\lim_{y \rightarrow + \infty} \frac{f(y)}{y^n}$ exists''. Remark that this new condition is not automatically verified for $n = 0$ (contrary to the first one).
\item We can easily show that the function $\tau$ involved in Theorem \ref{t4} is $o(x^n)$ (as $x \rightarrow + \infty$). For $n = 0$, this information implies that $\tau$ is a zero function (this is why $\tau$ doesn't appear in Theorem \ref{t3}); but for $n \geq 1$, the fact $\tau(x) = o(x^n)$ is obvious (since $\tau$ is continuous and $1$-periodic, so bounded).
\end{enumerate}
\end{rmks}

The statement of Theorem \ref{t4} above incites us to solve Equation $(\E_{n , a})$ ($n , a \in \N_0$, $a \geq 2$) in the set of real $1$-periodic functions. In what follows, we propose more precisely to solve that equation in the set of real $1$-periodic functions which are expandable as Fourier series. In this direction, we have obtained the following result:

\begin{thm}\label{t5}
Let $n \in \N$ and $a \geq 2$ be an integer. Let also $f$ be a real $1$-periodic function which is expandable as a Fourier series. Then $f$ is a solution of Equation $(\E_{n , a})$ if and only if it has the form{\rm:}
$$
f(x) = \sum_{k = 1}^{+ \infty} \left(u_k \frac{\cos(2 \pi k x)}{k^n} + v_k \frac{\sin(2 \pi k x)}{k^n}\right) ,
$$
where ${(u_k)}_{k \geq 1}$ and ${(v_k)}_{k \geq 1}$ are real sequences which satisfy the relations{\rm:}
$$
\left\{\begin{array}{rcl}
u_{a k} & = & u_k , \\
v_{a k} & = & v_k
\end{array}
\right. ~~~~~~~~~~ (\forall k \in \N) .
$$
\end{thm}

\begin{proof}
For simplicity, we use the exponential form of Fourier series. Let
$$
f(x) = \sum_{k \in \Z} \alpha_k e^{2 \pi i k x}
$$
be the Fourier series expansion of $f$, where ${(\alpha_k)}_{k \in \Z}$ is a sequence of complex numbers. Using the elementary formula
$$
e^{2 \pi i k x} + e^{2 \pi i k \left(x + \frac{1}{a}\right)} + e^{2 \pi i k \left(x + \frac{2}{a}\right)} + \dots + e^{2 \pi i k \left(x + \frac{a - 1}{a}\right)} = \begin{cases}
a e^{2 \pi i k x} & \text{if } k \equiv \restmod{0}{a} \\
0 & \text{else} 
\end{cases}
$$
(valid for all $k \in \Z$ and all $x \in \R$), we find that we have for all $x \in \R$:
\begin{align*}
f(x) + f\left(x + \frac{1}{a}\right) + f\left(x + \frac{2}{a}\right) + \dots + f\left(x + \frac{a - 1}{a}\right) & = \sum_{\begin{subarray}{c}
k \in \Z \\
k \equiv \restmod{0}{a}
\end{subarray}} \alpha_k a e^{2 \pi i k x} \\
& = \sum_{\ell \in \Z} a \alpha_{a \ell} e^{2 \pi i a \ell x} ~ (\text{by putting } k = a \ell) .
\end{align*}
Thus (for all $x \in \R$):
\begin{equation}\label{eq15}
a^{n - 1} \left(f(x) + f\left(x + \frac{1}{a}\right) + f\left(x + \frac{2}{a}\right) + \dots + f\left(x + \frac{a - 1}{a}\right)\right) = \sum_{\ell \in \Z} a^n \alpha_{a \ell} e^{2 \pi i a \ell x} .
\end{equation}
On the other hand, we have for all $x \in \R$:
\begin{equation}\label{eq16}
f(a x) = \sum_{k \in \Z} \alpha_k e^{2 \pi i a k x} .
\end{equation}
By comparing between \eqref{eq15} and \eqref{eq16}, we derive that:
$$
f \text{ is a solution of } (\E_{n , a}) \Longleftrightarrow \alpha_{a k} = \frac{\alpha_k}{a^n} ~~~~ (\forall k \in \Z) .
$$
For $k = 0$, the relation $\alpha_{a k} = \frac{\alpha_k}{a^n}$ gives $\alpha_0 = 0$. Next, for $k \in \Z^*$, by setting $\alpha_k = \frac{\beta_k}{k^n}$ (so ${(\beta_k)}_{k \in \Z^*}$ is a sequence of complex numbers), we derive that:
$$
f \text{ is a solution of } (\E_{n , a}) \Longleftrightarrow \alpha_0 = 0 \text{ and } \beta_{a k} = \beta_k ~ (\forall k \in \Z^*) .
$$
Finally, by setting for all $k \in \N$:
$$
u_k := 2 \, \Re(\beta_k) ~~\text{and}~~ v_k := - 2 \, \Im(\beta_k)
$$
and taking into account that $f$ is real, the Fourier series expansion of $f$ becomes
\begin{align*}
f(x) & = \sum_{k \in \Z} \alpha_k e^{2 \pi i k x} \\
& = \alpha_0 + \sum_{k \in \Z^*} \beta_k \frac{e^{2 \pi i k x}}{k^n} \\
& = \alpha_0 + \sum_{k = 1}^{+ \infty} \left(u_k \frac{\cos(2 \pi k x)}{k^n} + v_k \frac{\sin(2 \pi k x)}{k^n}\right)
\end{align*}
and we have
$$
f \text{ is a solution of } (\E_{n , a}) \Longleftrightarrow \alpha_0 = 0 , u_{a k} = u_k ~ (\forall k \in \N), \text{ and } v_{a k} = v_k ~ (\forall k \in \N) . 
$$
The theorem is proved.
\end{proof}

\begin{expls}
Let $n$ be a fixed positive integer which we take $\geq 2$ if necessary. Leaning on Theorem \ref{t5}, let us give some examples of real functions which are solutions of some Equations $(\E_{n , a})$ without being solutions of others. Denoting by $s_2(k)$ the sum of the binary digits of a given positive integer $k$, it is clear that we have $s_2(2 k) = s_2(k)$ for every $k \in \N$. So, according to Theorem \ref{t5}, the real function defined by:
$$
f_1(x) := \sum_{k = 1}^{+ \infty} \frac{s_2(k)}{k^n} \cos(2 \pi k x) ~~~~~~~~~~ (\forall x \in \R)
$$
is a solution of Equation $(\E_{n , 2})$ but it is not a solution of Equation $(\E_{n , 3})$ for example (since $s_2(3 k) \neq s_2(k)$ in general). Next, denoting by $\odd(k)$ the odd part of a given positive integer $k$ (i.e., the greatest odd positive divisor of $k$), the real function defined by:
$$
f_2(x) := \sum_{k = 1}^{+ \infty} \frac{\odd(k)}{k^n} \cos(2 \pi k x) ~~~~~~~~~~ (\forall x \in \R) 
$$
is a solution of Equation $(\E_{n , 2})$ (since we have obviously $\odd(2 k) = \odd(k)$ for any $k \in \N$), but it is not a solution of any other Equation $(\E_{n , a})$ if $a$ is not a power of $2$. Finally, denoting by ${(\mu_k)}_{k \in \N}$ the sequence defined by:
$$
\mu_k := \begin{cases}
1 & \text{if } k \text{ is a power of } 2 \\
0 & \text{else}
\end{cases} ~~~~~~~~~~ (\forall k \in \N)
$$
(which clearly satisfies the condition $\mu_{2 k} = \mu_k$, for all $k \in \N$), we obtain that the real function given by:
$$
f_3(x) := \sum_{k = 1}^{+ \infty} \frac{\mu_k}{k^n} \cos(2 \pi k x) = \sum_{\ell = 0}^{+ \infty} \frac{\cos\left(2^{\ell + 1} \pi x\right)}{2^{\ell n}} ~~~~~~~~~~ (\forall x \in \R)
$$
is a solution of Equation $(\E_{n , 2})$ but it is not a solution of any other Equation $(\E_{n , a})$ for $a$ not a power of $2$. Remark also that $f_3 \in \mathscr{C}^{n - 1}(\R , \R)$ but $f_3 \not\in \mathscr{C}^n(\R , \R)$. 
\end{expls}

Now, for a given $n \in \N$, we are going to characterize, among the real $1$-periodic functions expandable as Fourier series, those which are solutions of all Equations $(\E_{n , a})$ ($a \in \N$). We have the following corollary:

\begin{coll}\label{coll2}
Let $n \in \N$ and $f$ be a real $1$-periodic function which is expandable as a Fourier series. Then $f$ is a solution of all Equations $(\E_{n , a})$ {\rm(}$a \geq 2$ an integer{\rm)} if and only if it has the form{\rm:}
\begin{equation}\label{eq17}
f(x) = \alpha \sum_{k = 1}^{+ \infty} \frac{\cos(2 \pi k x)}{k^n} + \beta \sum_{k = 1}^{+ \infty} \frac{\sin(2 \pi k x)}{k^n}
\end{equation}
{\rm(}where $\alpha , \beta \in \R${\rm)}, or equivalently the form{\rm:}
$$
f(x) = \lambda B_n(x) + \mu \overline{B}_n(x) ~~~~~~ (\text{valid for all $x \in [0 , 1)$ if $n \geq 2$ and all $x \in (0 , 1)$ if $n = 1$})
$$
{\rm(}where $\lambda , \mu \in \R${\rm)}.
\end{coll}

\begin{proof}
By Theorem \ref{t5}, it is clear that the functions of the form in \eqref{eq17} are solutions of all Equations $(\E_{n , a})$ ($a \geq 2$). Conversely, suppose that $f$ is a solution of all Equations $(\E_{n , a})$ ($a \geq 2$) and let us show that $f$ has the form in \eqref{eq17}. By Theorem \ref{t5}, $f$ would have the form:
$$
f(x) = \sum_{k = 1}^{+ \infty} \left(u_k \frac{\cos(2 \pi k x)}{k^n} + v_k \frac{\sin(2 \pi k x)}{k^n}\right) ,
$$
where ${(u_k)}_{k \geq 1}$ and ${(v_k)}_{k \geq 1}$ are real sequences satisfying the relations
$$
\left\{\begin{array}{rcl}
u_{a k} & = & u_k \\
v_{a k} & = & v_k
\end{array}
\right. ~~~~~~~~~~ (\forall k , a \in \N) .
$$
In particular (by specializing $k$ to $1$), we find that for any $a \in \N$: $u_a = u_1$ and $v_a = v_1$. In other words, the sequences ${(u_k)}_k$ and ${(v_k)}_k$ are constant. Consequently, we have
$$
f(x) = u_1 \sum_{k = 1}^{+ \infty} \frac{\cos(2 \pi k x)}{k^n} + v_1 \sum_{k = 1}^{+ \infty} \frac{\sin(2 \pi k x)}{k^n} ,
$$
which is the required form of $f$. The last part of the corollary is an immediate consequence of Formulas \eqref{eq18} and \eqref{eq19} and some elementary trigonometric formulas. This completes the proof of the corollary.
\end{proof}

\begin{expl}
The Fourier series expansion of the real $1$-periodic function $x \mapsto \log\left(2 \vert\sin(\pi x)\vert\right)$ is well known and given by:
$$
\log\left(2 \left\vert\sin(\pi x)\right\vert\right) = - \sum_{k = 1}^{+ \infty} \frac{\cos(2 \pi k x)}{k} ~~~~~~~~~~ (\forall x \in \R \setminus \Z) .
$$
So, according to Corollary \ref{coll2}, this ``elementary'' function is a solution of all Equations $(\E_{1 , a})$ ($a \in \N$).
\end{expl}

\rhead{\textcolor{OrangeRed3}{\it References}}

\end{document}